\documentstyle[fleqn,12pt]{article}
\begin{document}
\pagenumbering{arabic}\setcounter{page}{1}
\pagestyle{plain}\baselineskip=16pt

\thispagestyle{empty}
\vspace{1.4cm}

\begin{center}
{\Large\bf Differential Geometry of the Lie algebra of the
quantum plane}
\end{center}

\vspace{1cm}
\begin{center} Salih \c Celik \footnote{E-mail address: sacelik@yildiz.edu.tr}
 and Sultan A \c Celik

Yildiz Technical University, Department of Mathematics, \\
34210 Davutpasa-Esenler, Istanbul, TURKEY. \end{center}

\vspace{2cm} {\small We present a differential calculus on the
extension of the quantum plane obtained considering that the
(bosonic) generator $x$ is invertible and furthermore working
polynomials in $\ln x$ instead of polynomials in $x$. We call
quantum Lie algebra to this extension and we obtain its Hopf
algebra structure and its dual Hopf algebra.}

\vfill\eject\noindent

Quantum plane [1] is a simple example of quantum space and has been studied
intensively by many authors in the past years. It can be obtained by
deformation of the classical plane [2]. For references to the literature we refer
to the recent book by Majid [3].

As usual in noncommutative geometry [4] quantum plane has many differential
calculi. By using the noncommutative differential geometry method in Ref. 5
and interpreting the dual plane of the quantum plane to consist of
differentials of the coordinates of the quantum plane, covariant differential
calculus on the quantum plane has been developed. Therefore the quantum
plane provides a simple example for noncommutative differential geometry [6].

It is known that, in order to construct a noncommutative differential calculus
on quantum groups and Hopf algebras, one takes into consideration the
associative algebra of functions on the group. The starting point of this work
is its Lie algebra. For quamtum superplane, similar a work studied in Ref. 7.
In this work, we present here a differential calculus on the Lie algebra
of the associative algebra of functions on the extended quantum plane.



The quantum plane [1] is defined as an associative algebra generated by two
noncommuting coordinates $x$ and $y$ with the single quadratic relation
$$ x y - q y x = 0 \qquad q \in {\cal C} - \{0\}. \eqno(1) $$
This algebra is known as the algebra of polynomials over the
quantum plane and we shall denote by ${\cal A}_q = Fun_q(R(2))$.
We extend the algebra ${\cal A}_q$ by including inverse of $x$
which obeys
$$x x^{-1} = 1 = x^{-1} x.$$

We know, from Ref. 8, that the extended algebra ${\cal A}_q$ is a
Hopf algebra with the following costructures:
$$\Delta(x) = x \otimes x \qquad
  \Delta(y) = y \otimes 1 + x \otimes y $$
$$\epsilon(x) = 1 \qquad \epsilon(y) = 0 \eqno(2)$$
$$\kappa(x) = x^{-1} \qquad \kappa(y) = - x^{-1} y. $$
It is not difficult to verify the following properties of
costructures:
$$(\Delta \otimes \mbox{id}) \circ \Delta =
  (\mbox{id} \otimes \Delta) \circ \Delta $$
$$m \circ (\epsilon \otimes \mbox{id}) \circ \Delta
  = m \circ (\mbox{id} \otimes \epsilon) \circ \Delta \eqno(3)$$
$$m \circ (\kappa \otimes \mbox{id}) \circ \Delta = \epsilon
  = m \circ (\mbox{id} \otimes \kappa) \circ \Delta $$
where id denotes the identity mapping and $m$ is the multiplication
map $$m : {\cal A}_q \otimes {\cal A}_q \longrightarrow {\cal A}_q \qquad
  m(a \otimes b) = ab. $$


It is known that an element of a Lie group can be represented by exponential
of an element of its Lie algebra. Using this fact, we can define the
generators of ${\cal A}_q$ as [9]
$$x = e^X \qquad y = e^X Y \qquad q = e^h. \eqno(4)$$
Then we obtain the relation
$$[X,Y] = h Y \eqno(5)$$
where
$$[a,b]_\pm = a b \pm ba. $$
This is the relations of a Lie algebra and we shall denote it by
${\cal L(A}_q) =: g_h$. The Hopf algebra structure of $g_h$ can be
read off from (2),
$$\Delta(X) = X \otimes 1 + 1 \otimes X $$
$$\Delta(Y) = Y \otimes e^{-X} + 1 \otimes Y $$
$$\epsilon(X) = 0 \qquad \epsilon(Y) = 0 \eqno(6)$$
$$\kappa(X) = - X \qquad \kappa(Y) = - Y e^X. $$


We now shall build up a noncommutative differential calculus on the Lie
algebra $g_h$. This must be involve functions on the Lie algebra $g_h$,
differentials and differential forms. For this, we have to define a linear
map $\delta$ which acts on the functions of the elements of $g_h$.
The exterior differential $\delta$ is an operator which gives the mapping
from the generators of $g_h$ to the differentials:
$$\delta : u \mapsto \delta u \qquad u \in \{X,Y\}.$$
We demand that the exterior differential $\delta$ has to satisfy two
properties: the nilpotency
$$\delta^2 = 0 \eqno(7a)$$
and the Leibniz rule
$$\delta(f g) = (\delta f) g + f (\delta g). \eqno(7b)$$

In order to establish a noncommutative differential calculus on the algebra
$g_h$, we assume that the commutation relations between the elements of $g_h$
and their differentials are the following form:
$$X ~\delta X = A_{11} \delta X ~X + B_1 \delta X + B_2 \delta Y $$
$$X ~\delta Y = A_{12} \delta Y ~X + B_3 \delta X + B_4 \delta Y $$
$$Y ~\delta X = A_{21} \delta X ~Y + B_5 \delta X + B_6 \delta Y $$
$$Y ~\delta Y = A_{22} \delta Y ~Y + B_7 \delta X + B_8 \delta Y. \eqno(8)$$
The coefficients $A_{ij}$ and $B_i$ will be determined in terms of the
"new" deformation parameter $h$. To find them we shall use the
consistency of calculus.

Using the consistency conditions
$$\delta (X Y - Y X - h Y) = 0 $$
$$(X Y - Y X - h Y) \delta X = 0 $$
$$(X Y - Y X - h Y) \delta Y = 0 $$
we obtain the system 
$$A_{11} = A_{12} = A_{21} = 0 \qquad B_2(1 - A_{22}) = 0$$
$$B_3 - B_5 = 0 \qquad B_4 - B_6 = h \qquad B_3 B_6 - B_2 B_7 = h B_5$$
$$(B_1 - B_4) B_6 + (B_8 - B_5) B_2 = - h B_6$$
$$B_3 B_6 - B_2 B_7 = h B_8$$
$$(B_1 - B_4) B_7 + (B_8 - B_5) B_3 =  h B_7$$
or 
$$A_{11} = A_{12} = A_{21} = 0 \qquad B_2(1 - A_{22}) = 0$$
$$B_3 = B_5 = B_8 \qquad B_4 - B_6 = h \qquad (B_1 - B_4 - h) B_7 = 0.$$
Thus we have the following algebra 
$$[X, \delta X] = B_1 \delta X + B_2 \delta Y $$
$$[X, \delta Y] = B_3 \delta X + B_4 \delta Y $$
$$[Y, \delta X] = B_3 \delta X + (B_4 - h) \delta Y $$
$$[Y, \delta Y] = (A_{22} - 1) \delta Y ~Y + B_7 \delta X + B_3 \delta Y. 
\eqno(9)$$
But we do not have the required relations because of the constants. To find 
the constans in the above relations, we shall consider the covariance of the 
noncommutative differential calculus. 

We first note that consistency of a differential calculus with commutation 
relation (5) means that the algebra is a graded associative algebra 
generated by the elements of the set $\{X,Y, \delta X,\delta Y\}$. 
So, it is sufficient only describe the actions of co-maps on the subset 
$\{\delta X, \delta Y\}$. 

Let $\Gamma$ be a free left module over the algebra $g_h$ generated by the 
elements of the set $\{X,Y, \delta X, \delta Y\}$. The module $\Gamma$ 
becomes a unital associative algebra if one defines a multiplication law 
on $\Gamma$ by relations (5) and (9). 

We consider a map 
$\phi_L : \Gamma \longrightarrow g_h \otimes \Gamma$ such that 
$$\phi_L \circ \delta = (\mbox{id} \otimes \delta) \circ \Delta. \eqno(10)$$
Thus we have 
$$\phi_L (\delta X) = 1 \otimes \delta X \qquad 
  \phi_L (\delta Y) = 1 \otimes \delta Y + 
  {{1 - e^h}\over h} ~Y \otimes X e^{-X}. \eqno(11)$$
We now define a map $\Delta_L$ as follows: 
$$\Delta_L(a_1 ~\delta b_1 + \delta b_2 ~a_2) = 
  \Delta(a_1) \phi_L(\delta b_1) + \phi_L(\delta b_2) 
  \Delta(a_2). \eqno(12)$$
Applying the linear map $\Delta_L$ to relations (9) we get, as the solution 
of the obtained system, 
$$A_{11} = A_{12} = A_{21} = 1 \qquad A_{22} = e^h$$
$$B_1 = B_6 = - h \qquad B_2 = B_3 = B_4 = B_5 = B_7 = B_8 = 0. \eqno(13)$$
Consequently, the relations (9) is of the form 
$$[X, \delta X] = - h \delta X \qquad [X, \delta Y] = 0$$
$$[Y, \delta X] = - h \delta Y \qquad 
  [Y, \delta Y] = (e^h - 1) \delta Y ~Y. \eqno(14)$$
Applying the exterior differential $\delta$ to the relation (14) one has 
$$[\delta X, \delta Y]_+ = 0 \qquad (\delta X)^2 =0= (\delta Y)^2. \eqno(15)$$
Note that, the differentials of $X$ and $Y$ in terms of $x$ and $y$ as
follows:
$$\delta X = {h\over {1 - e^{-h}}} ~\delta x ~x^{-1} \qquad
  \delta Y = x^{-1} (\delta y - \delta x ~x^{-1} y).$$

We now consider another map $\phi_R: \Gamma \longrightarrow \Gamma \otimes g_h$ 
such that 
$$\phi_R \circ \delta = (\delta \otimes \mbox{id}) \circ \Delta. \eqno(16)$$
Thus one has 
$$\phi_R (\delta X) = \delta X \otimes 1 \qquad 
  \phi_R (\delta Y) = \delta Y \otimes e^{-X}. \eqno(17)$$
We define a map $\Delta_R$ with again (12) by replacing $R$ with $L$. The 
map $\Delta_R$ also leaves invariant the relations (14) and (15), and the 
following identities are satisfied: 
$$(\Delta \otimes \mbox{id}) \circ \Delta_L = 
  (\mbox{id} \otimes \Delta_L) \circ \Delta_L 
  \qquad (\epsilon \otimes \mbox{id}) \circ \Delta_L = \mbox{id} $$
$$(\mbox{id} \otimes \Delta) \circ \Delta_R = 
  (\Delta_R \otimes \mbox{id}) \circ \Delta_R 
  \qquad (\mbox{id} \otimes \epsilon) \circ \Delta_R = \mbox{id}. \eqno(18)$$

To denote the coproduct, counit and coinverse which will be defined on the 
algebra $\Gamma$ with those of ${\cal A}$ may be inadvisable. For this reason, 
we shall denote them with a different notation. Let us define the 
map $\hat{\Delta}$ as 
$$\hat{\Delta} = \Delta_R + \Delta_L \eqno(19)$$
which will allow us to define the coproduct of the differential algebra. 
We denote the restriction of $\hat{\Delta}$ to the algebra $g_h$ by 
$\Delta$ and the extension of $\Delta$ to the differential algebra $\Gamma$ 
by $\hat{\Delta}$. It is possible to interpret the relation 
$$\hat{\Delta}\vert_{g_h} = \Delta \eqno(20)$$
as the definition of $\hat{\Delta}$ on the generators of $g_h$ and 
(19) as the definition of $\hat{\Delta}$ on differentials. 

It is not difficult to verify the following conditions: 

{\bf a)} $\Gamma$ is an $g_h$-bimodule, 

{\bf b)} $\Gamma$ is an $g_h$-bicomodule with left and right coactions 
$\Delta_R$ and $\Delta_L$, respectively, making $\Gamma$ a left 
and right $g_h$-comodule with (18), and 
$$(\Delta_L \otimes \mbox{id}) \circ \Delta_R = 
  (\mbox{id} \otimes \Delta_R) \circ \Delta_L \eqno(21)$$
which is the $g_h$-bimodule property. So, the triple 
$(\Gamma, \Delta_L, \Delta_R)$ is a bicovariant bimodule over Hopf algebra 
$g_h$. In additional, since 

{\bf c)} $(\Gamma, \delta)$ is a first order differential calculus over 
$g_h$, and 

{\bf d)} $\delta$ is both a left and a right comodule map, i.e. for all 
$a \in g_h$ 
$$(\mbox{id} \otimes \delta) ~\Delta(a) = \Delta_L(\delta a) \qquad 
  (\delta \otimes \mbox{id}) ~\Delta(a) = \Delta_R(\delta a) 
  \eqno(22)$$
the quadruple 
$(\Gamma, d, \Delta_L, \Delta_R)$ is a first order bicovariant differential 
calculus over Hopf algebra $g_h$. 
 
Now let us return Hopf algebra structure of $\Gamma$. If we define a counit 
$\hat{\epsilon}$ for the differential algebra as 
$$\hat{\epsilon} \circ \delta = \delta \circ \epsilon = 0 \eqno(23)$$
and 
$$\hat{\epsilon}\vert_{g_h} = \epsilon \qquad 
  \epsilon\vert_\Gamma = \hat{\epsilon}. \eqno(24)$$
we have 
$$\hat{\epsilon}(\delta X) = 0 \qquad \hat{\epsilon}(\delta Y) = 0 
  \eqno(25)$$
where 
$$\hat{\epsilon}(a_1 ~\delta b_1 + \delta b_2 ~a_2) = 
  \epsilon(a_1) \hat{\epsilon}(\delta b_1) + 
  \hat{\epsilon}(\delta b_2) \epsilon(a_2). \eqno(26)$$
Here we used the fact that $\delta (1) = 0$. 

As the next step we obtain a coinverse $\hat{\kappa}$. For this, 
it suffices to define $\hat{\kappa}$ such that 
$$\hat{\kappa} \circ \delta = \delta \circ \kappa \eqno(27)$$
and 
$$\hat{\kappa}\vert_{g_h} = \kappa \qquad 
  \kappa\vert_\Gamma = \hat{\kappa} \eqno(28)$$
where 
$$\hat{\kappa}(a_1 ~\delta b_1 + \delta b_2 ~a_2) = 
  \hat{\kappa}(\delta b_1) \kappa(a_1) + \kappa(a_2) \hat{\kappa}(\delta b_2). 
\eqno(29)$$
Thus the action of $\hat{\kappa}$ on the generators $\delta X$ and 
$\delta Y$ is as follows: 
$$\hat{\kappa}(\delta X) = - \delta X \qquad 
 \hat{\kappa}(\delta Y) = e^{X-h} \left(\delta Y + 
  {{e^h - 1}\over h} ~\delta X ~Y\right).   \eqno(30)$$
Note that it is easy to check that $\hat{\epsilon}$ and $\hat{\kappa}$ leave 
invariant the relations (14) and (15). 
Consequently, we can say that the structure 
$(\Gamma, \hat{\Delta}, \hat{\epsilon}, \hat{\kappa})$ 
is a graded Hopf algebra. 


Just as we introduced the derivatives of the
generators of ${\cal A}_q$ in the standard way, let us introduce derivatives
of the generators of $g_h$ and multiply explicit expression of the exterior
differential $\delta$ from the right by $X f$ and $Y f$, respectively. Then,
using the Leibniz rule for partial derivatives
$$\partial_i (f g) = (\partial_i f) g + f (\partial_i g) $$
we obtain
$$[\partial_X, X] = {h\over {1 - e^{-h}}} \left(1 + 
   (e^{-h} - 1) \partial_X \right) $$
$$[\partial_X, Y] = 0 \qquad [\partial_Y, X] = 0 $$
$$[\partial_Y, Y] = e^h + (1 - e^h) \left(\partial_X - Y \partial_Y \right).
  \eqno(31)$$
Note that, the partial derivatives of $X$ and $Y$ in terms of $x$
and $y$ as follows:
$$\partial_X = x \partial_x + y \partial_y \qquad
\partial_Y = e^h x \partial_y.$$
The commutation relations among the partial derivatives can be
easily obtained by using $\delta^2 = 0$. So it follows that
\begin{eqnarray*}
 \delta & = & {{1 - e^{-h}}\over h} ~ \delta X ~\partial_X + e^{-h}
            \delta Y ~\partial_Y\\
 0 & = & \delta^2 = \delta \left({{1 - e^{-h}}\over
h} ~\delta X ~\partial_X +
       e^{-h} \delta Y \partial_Y \right) \\
& = & e^{-h} ~{{e^{-h} - 1}\over h} ~\delta X \delta Y
      \left(\partial_Y \partial_X - \partial_X \partial_Y\right)
\end{eqnarray*}
which says that
$$[\partial_X, \partial_Y] = 0. \eqno(32)$$

In order to find the commutation relations between the differentials and
partial derivatives, we demand a form similar to (8). After some tedious but
straigtword calculations, we find
$$[\partial_X, \delta X] = (e^h - 1) \delta X \partial_X -
   e^h \delta X $$
$$[\partial_X, \delta Y] = (e^h - 1) \delta Y \partial_X -
   e^h \delta Y $$
$$[\partial_Y, \delta X] = 0 \eqno(33)$$
$$[\partial_Y, \delta Y] = (e^h - 1) \delta Y \partial_Y +
{{(e^h - 1)^2}\over h} ~\delta X \partial_X + e^h ~{{1-e^h}\over h} \delta X.$$


We now define two one-forms using the generators of $g_h$
and show that the algebra of forms is a graded Hopf algebra.

If we call them $w_1$ and $w_2$ then we can define them as follows:
$$w_1 = {{1 - e^{-h}}\over h} ~\delta X \qquad w_2 = e^X \delta Y. \eqno(34)$$
We denote the algebra of forms generated by $w_1$ and $w_2$ by $\Omega$.
The generators of the algebra $\Omega$ with the generators of $g_h$ satisfy
the following relations:
$$[w_1, X] = h w_1 \qquad [w_1, Y] = (1 - e^{-h}) e^{-X} w_2 $$
$$[w_2, X] = 0 \qquad [w_2, Y] = 0. \eqno(35)$$
The commutation rules of the generators of $\Omega$ are
$$w_1 w_2 = - e^h w_2 w_1 \qquad w_1^2 = 0 = w_2^2. \eqno(36)$$

One can make the algebra $\Omega$ into a graded Hopf algebra with the
following costructures: the coproduct
$\Delta: \Omega \longrightarrow \Omega \otimes \Omega$ is given by
$$\Delta(w_1) = w_1 \otimes 1 + 1 \otimes w_1 $$
$$\Delta(w_2) = w_2 \otimes 1 + e^X \otimes w_2 - e^X Y \otimes w_1.\eqno(37)$$
The counit $\epsilon: \Omega \longrightarrow {\cal C}$ is given by
$$\epsilon(w_1) = 0 \qquad \epsilon(w_2) = 0 \eqno(38)$$
and the coinverse $\kappa: \Omega \longrightarrow \Omega$ is defined by
$$\kappa(w_1) = - w_1 \qquad \kappa(w_2) = e^{-X} w_2 + Y w_1. \eqno(39)$$


Now we shall obtain the universal enveloping algebra of
$g_h$, ${\cal U}(g_h)$, via the vector fields related to the one-forms.
We first write the Cartan-Maurer forms as
$$\delta X = {h\over {1-e^{-h}}} ~w_1 \qquad \delta Y = e^{-X} w_2. \eqno(40)$$
Writing the exterior differential $\delta$ in the form
$$\delta = w_1 T_1 + w_2 T_2$$
and considering an arbitrary function $f$ of the generators of $g_h$
and using the nilpotency of $\delta$ one has
$$T_2 T_1 = e^{-h} T_1 T_2 + T_2. \eqno(41)$$
The commutation relations of vector fields with the generators of
$g_h$ as follows:
$$[T_1, X] = {h\over {1 - e^{-h}}} \left(1 + (e^{-h} - 1) T_1\right) $$
$$[T_1, Y] = 0 \qquad [T_2, X] = 0 $$
$$[T_2, Y] = e^{-X} \left(1 + (e^{-h} - 1) T_1\right). \eqno(42)$$

The Hopf algebra structure of the vector fields is given by
$$\Delta(T_1) = T_1 \otimes 1 + 1 \otimes T_1 + (e^{-h} - 1) T_1 \otimes T_1 $$
$$\Delta(T_2) = T_2 \otimes 1 + 1 \otimes T_2 + (e^{-h} - 1) T_2 \otimes T_1 $$
$$\epsilon(T_1) = 0 \qquad \epsilon(T_2) = 0 \eqno(43)$$
$$\kappa(T_1) = - T_1 \left(1 + (e^{-h} - 1) T_1 \right)^{-1} $$
$$\kappa(T_2) = - T_2 \left\{1 - (e^{-h} - 1) T_1 \left(1 + (e^{-h} - 1) T_1\right)^{-1}\right\}.$$

If we introduce the operators $H$ and $N$ as
$$T_1 = {{e^{-h N} - 1}\over {e^{-h} - 1}} \qquad \mbox{and} \qquad
   T_2 = e^{-h N} H \eqno(44)$$
we get
$$[H, N] = H. \eqno(45)$$
Then, dual Hopf algebra to ${\cal U}(g_h)$ is given by
$$\Delta(N) = N \otimes 1 + 1 \otimes N $$
$$\Delta(H) = H \otimes 1 + e^{h N} \otimes H $$
$$\epsilon(N) = 0 \qquad \epsilon(H) = 0 \eqno(46)$$
$$\kappa(N) = - N \qquad \kappa(H) = - e^{-h N} H.$$

\noindent
{\bf Acknowledgment}

\noindent
This work was supported in part by TBTAK the Turkish Scientific and Technical
Research Council.


\end{document}